\pdfoutput=1  
\documentclass {amsart}  

\usepackage{amsmath}

\usepackage{amsfonts}

\usepackage{amssymb}
\usepackage[bookmarks=true, colorlinks=true, linkcolor=black, citecolor=black, filecolor=black, urlcolor=black, pdfstartview=FitH]{hyperref}

\usepackage[ansinew]{inputenc}

\newtheorem{theorem}{Theorem}
\newtheorem{lemma}{Lemma}
\newtheorem{corollary}{Corollary}
\newtheorem{proposition}{Proposition}

\def\qmod#1{\ ({\rm mod}\ #1)}

\def\ord{\rm ord}

\begin{document}

\title{Strengthening the Baillie-PSW primality test}

\author{Robert Baillie}
\address{State College, PA 16803-3029 USA}
\email[Robert Baillie]{rjbaillie@frii.com}
\author{Andrew Fiori}
\address{Mathematics and Computer Science \\
  4401 University Drive \\
  University of Lethbridge \\
  Lethbridge, Alberta, T1K 3M4 Canada}
\email[Andrew Fiori]{andrew.fiori@uleth.ca}
\thanks{A.F.'s work was supported partially by the University of Lethbridge and NSERC}
\author{Samuel S. Wagstaff, Jr.}
\address{Center for Education and Research in Information Assurance and Security
  and Department of Computer Sciences, Purdue University \\
  West Lafayette, IN 47907-1398 USA}
\email[Samuel S. Wagstaff, Jr.]{ssw@cerias.purdue.edu}
\thanks{S.S.W.'s work was supported by the CERIAS Center at Purdue University}

\subjclass[2010]{Primary 11Y11; Secondary 11A51} 
\keywords{primality test, Lucas sequences}

\begin{abstract}
In 1980, the first and third authors proposed a probabilistic primality test that has become known as the Baillie-PSW (BPSW) primality test.
Its power to distinguish between primes and composites comes from combining a Fermat probable prime test with a Lucas probable prime test.
No odd composite integers have been reported to pass this combination of primality tests if the parameters are chosen in an appropriate way.
Here, we describe a significant strengthening of this test that comes at almost no additional computational cost.
This is achieved by including in the test Lucas-V pseudoprimes, of which there are only five less than $10^{15}$.
\end{abstract}

\maketitle

\section{Introduction}

A (Fermat) \emph{base-}$a$ \emph{pseudoprime}, or psp($a$), is a composite positive integer $n$ that satisfies the conclusion of Fermat's little theorem, that is
\[
a^{n-1}\equiv 1 \qmod{n} \thinspace .
\]
For each integer base $a > 1$, there are infinitely many pseudoprimes, but they are sparser than primes.
In \cite{BW80}, the first and third authors studied analogues of pseudoprimes in which $a^{n-1}-1$ is replaced by a Lucas sequence.

Let $D$, $P$ and $Q$ be integers with $P>0$ and $D = P^2 - 4Q \neq 0$.
Define $U_0 = 0$, $U_1 = 1$, $V_0 = 2$ and $V_1 = P$.
The Lucas sequences $U_k$ and $V_k$ with parameters $P$ and $Q$ are defined for $k\ge2$ by
\[
U_k = P U_{k-1} - Q U_{k-2} \quad \text{and} \quad V_k = P V_{k-1} - Q V_{k-2} \thinspace .
\]

Let $n > 1$ be an odd positive integer.
Choose $D$, $P$, and $Q$ so that the Jacobi symbol $(D/n) = -1$.
It is well known \cite{BW80}, \cite{BLS1975} that if $n$ is prime and $(n, Q) = 1$, then
\begin{align}
U_{n+1} & \equiv 0 \qmod{n} \thinspace , \label{E:UCongruence} \\
V_{n+1} & \equiv 2Q \qmod{n} \thinspace . \label{E:VCongruence}
\end{align}

In \cite{BW80}, we defined a \emph{Lucas pseudoprime} with parameters $P$ and $Q$ to be a \emph{composite} integer $n$ satisfying \eqref{E:UCongruence}.
We proposed a fast probable prime test by combining the Lucas primality criterion in \eqref{E:UCongruence} with a (Fermat) probable prime test.

In this paper, we emphasize the importance of the primality criterion in Congruence \eqref{E:VCongruence}.
We found that, using a standard method of choosing $D$, $P$, and $Q$,
among the composite $n$ under $10^{15}$, there are over two million that satisfy \eqref{E:UCongruence},
but only \emph{five} that satisfy \eqref{E:VCongruence}.

\newpage  

\textbf{Layout of this paper.}
\begin{itemize}
 \item Section \ref{S:Background}: we give details on Fermat and Lucas pseudoprimes and describe how to efficiently compute terms in the Lucas sequences;
 \item Section \ref{S:OriginalBPSW}: we define the original Baillie-PSW primality test, we list applications that use this test, and we summarize calculations that have been performed over the past 40 years;
 \item Section \ref{S:data15}: we summarize the data on pseudoprimes up to $10^{15}$;
 \item Section \ref{S:IsAStarSpecial}: we discuss whether the scarcity of composite solutions to \eqref{E:VCongruence} is due to the particular method for choosing $P$ and $Q$;
 \item Section \ref{S:newTest}: we propose a strengthened primality test that includes Congruence \eqref{E:VCongruence} and offer a reward for a counterexample;
 \item Section \ref{S:Qgt1}: discusses the importance of choosing $Q$ to be neither $+1$ nor $-1 \qmod{n}$;
 \item Section \ref{S:counterexamples}: we reprise Pomerance's heuristic argument that there are infinitely many counterexamples to the enhanced test;
 \item Appendix \ref{S:Astar}: we prove that two popular methods for choosing $P$ and $Q$ produce exactly the same Lucas pseudoprimes.
\end{itemize}

The authors thank Carl Pomerance for suggesting the proof of Theorem \ref{thmspspr} in Section \ref{S:counterexamples}.

\section{Background} \label{S:Background}

\subsection{Fermat probable primes and pseudoprimes} \label{S:IntroFermat}

A (Fermat) \emph{base-}$a$ \emph{probable prime}, or prp($a$), is a positive integer $n$ that satisfies the conclusion of Fermat's little theorem.
That is, if $n$ is prime and $(a, n) = 1$, then

\begin{equation} \label{E:FermatTest}
a^{n-1} \equiv 1 \qmod{n} \thinspace .
\end{equation}

The converse of Fermat's little theorem is not true, but if \eqref{E:FermatTest} is true for a given $a > 1$, then $n$ is likely to be prime.

A \emph{base-}$a$ \emph{pseudoprime}, or psp($a$), is a \emph{composite} $n$ that satisfies \eqref{E:FermatTest}.

Base-2 pseudoprimes up to $25 \cdot 10^9$ were studied in detail in \cite{PSW1980}.
The first ten base-2 pseudoprimes are 341, 561, 645, 1105, 1387, 1729, 1905, 2047, 2465, and 2701.

Since \cite{PSW1980} appeared in 1980, Feitsma \cite{Feitsma} has computed the psp(2) $< 2^{64} \approx 1.8 \cdot 10^{19}$.
There are 118\,968\,378 of them.

There are $\pi(2^{64}) - 1 = 425 \, 656 \, 284 \, 035 \, 217 \, 742$ odd primes $< 2^{64}$ \cite{Staple2015}.
Therefore, up to $2^{64}$, congruence \eqref{E:FermatTest} with $a = 2$ holds for $425 \, 656 \, 284 \, 035 \, 217 \, 742 + 118 \, 968 \, 378$ values of $n$, of which 99.9999999721 percent are prime.
This is why, if $2^{n-1} \equiv 1 \qmod{n}$, it is legitimate to call $n$ a \emph{probable prime}, and why this congruence is sometimes used as part of a test for primality.

Euler's criterion states that if $n$ is an odd prime and $(a, n) = 1$, then
\begin{equation}
a^{(n-1)/2} \equiv \left( \frac{a}{n} \right) \qmod{n} \thinspace , \label{E:EulerCriterion}
\end{equation}
where $\left( \frac{a}{n} \right)$ is the Jacobi symbol.
A composite number that satisfies this congruence is called a \emph{base-}$a$ \emph{Euler pseudoprime} (epsp($a$)).
The first ten epsp(2) are 561, 1105, 1729, 1905, 2047, 2465, 3277, 4033, 4681, and 6601.
The epsp($a$) are a proper subset of the psp($a$).
About half of the psp($a$) are epsp($a$) \cite[p. 1005]{PSW1980}, so \eqref{E:EulerCriterion} is a slightly stronger primality test than \eqref{E:FermatTest}.

\subsection{Strong probable primes and pseudoprimes} \label{S:IntroFermatStrong}

We now describe an even stronger, and more widely-used primality test, also based on Fermat's little theorem.
\cite{PSW1980} defines {\em strong} probable primes and {\em strong} pseudoprimes.
If $n$ is odd, then we can write $n - 1 = d \cdot 2^s$ where $d$ is odd.
If $n$ is an odd prime and $(a, n) = 1$, then either
\begin{align}
a^d & \equiv 1 \qmod{n}, \quad\text{or} \label{E:StrongPrp1} \\
a^{d \cdot 2^r} & \equiv -1 \qmod{n}, \quad\text{for some } r \text{ with } 0 \le r < s. \label{E:StrongPrp2}
\end{align}

If either \eqref{E:StrongPrp1} or \eqref{E:StrongPrp2} is true, then $n$ is called a \emph{base-}$a$ \emph{strong probable prime} (sprp($a$)).
If either of these holds, then we also have $a^{n-1} = a^{d \cdot 2^s} \equiv 1 \qmod{n}$.

If $n$ is composite and either \eqref{E:StrongPrp1} or \eqref{E:StrongPrp2} is true, then $n$ is called a \emph{base-}$a$ \emph{strong pseudoprime} (spsp($a$)).
The spsp($a$) are a proper subset of psp($a$), and so are scarcer than psp($a$).
For example, of the 118\,968\,378 psp(2) $< 2^{64}$ found by Feitsma, only 31\,894\,014 are spsp(2) \cite{Feitsma}.
The first ten base-2 strong pseudoprimes are 2047, 3277, 4033, 4681, 8321, 15841, 29341, 42799, 49141, and 52633.

The spsp($a$) are also a proper subset of epsp($a$).
Therefore, it makes sense for a primality test to use the strong conditions \eqref{E:StrongPrp1} and \eqref{E:StrongPrp2} instead of \eqref{E:FermatTest} or \eqref{E:EulerCriterion}.



To efficiently compute $a^{n-1}$, we use the binary expansion of $n - 1$.
The number of steps is essentially the number of \emph{binary digits} in $n$, that is, $\log_2(n)$.
We start by using the bits of $d$ to compute $a^d$:
Given the result $a^k$ from the previous step, we square $a^k$ to get $a^{2k}$, and if the corresponding bit is 1, we multiply by $a$ to get $a^{2k+1}$.
Once we reach $a^d$, the remaining $s$ bits of $n - 1$ are all 0, so we then square the result $s$ times.
We also reduce the result $\qmod{n}$ at each step to keep the sizes of the numbers reasonable.
This procedure fits nicely with the definition of a strong probable prime.
The following example shows how this works.

Example: Suppose $n = 341 = 11 \cdot 31$ and $a = 2$.
Then $n - 1 = 340 = 85 \cdot 2^2$, so $d = 85$ and $s = 2$.
In binary, $n - 1 = 101010100$ and $d = 1010101$.
The exponent, 340, has 9 bits, so there are (at most) 9 steps.
We process the bits from the left.

First, initialize the result $r$, to 1: $r \leftarrow 1$; this is $2^0$.
\begin{itemize}  
 \item bit 1 is 1: $r \leftarrow r^2 \cdot 2= 2^1 \equiv 2 \qmod{341}$ ;
 \item bit 2 is 0: $r \leftarrow r^2 = 2^2 \equiv 4 \qmod{341}$ ;
 \item bit 3 is 1: $r \leftarrow r^2 \cdot 2 = 2^5 \equiv 32 \qmod{341}$ ;
 \item bit 4 is 0: $r \leftarrow r^2 = 2^{10} \equiv 1 \qmod{341}$ ;
 \item bit 5 is 1: $r \leftarrow r^2 \cdot 2 = 2^{21} \equiv 2 \qmod{341}$ ;
 \item bit 6 is 0: $r \leftarrow r^2 = 2^{42} \equiv 4 \qmod{341}$ ;
 \item bit 7 is 1: $r \leftarrow r^2 \cdot 2 = 2^{85} \equiv 32 \qmod{341}$ ; this is $2^d$
 \item bit 8 is 0: $r \leftarrow r^2 = 2^{170} \equiv 1 \qmod{341}$ ;
 \item bit 9 is 0: $r \leftarrow r^2 = 2^{340} \equiv 1 \qmod{341}$ ;
\end{itemize}

The last step shows that 341 is a probable prime base 2.

However, note that, after using bit 7 to reach $2^{d} = 2^{85}$, we see that $2^d \not \equiv \pm1$.
In the next step, $2^{2d} = 2^{170} = 2^{(n-1)/2} \not \equiv -1 \qmod{341}$, so $n = 341$ is not a \emph{strong} probable prime base 2.
This proves that $n$ is composite, so the step with bit 9 is not needed for the strong test.
Therefore, the strong test, in addition to allowing fewer composites to pass, will also terminate at least one step earlier than if we simply computed $a^{n-1} \qmod{n}$.

A \emph{Carmichael number} is a composite integer $n$ that is a pseudoprime to every base $a$ for which $(a, n) = 1$.
They are also sparse, although there are infinitely many of them \cite{AGP-InfCarm1994}.
However, there are no \emph{strong} Carmichael numbers, that is, there is no composite $n$ which is \emph{strong} pseudoprime to all bases relatively prime to $n$:
Rabin proved \cite[Theorem 1]{Rabin1980} that any composite $n$ is a strong pseudoprime to at most $1/4$ of bases $a$, $1 \leq a < n$.

\subsection{Lucas sequences and pseudoprimes; Lucas-V pseudoprimes} \label{S:IntroLucas}

Lucas sequences, and their applications to prime-testing, were discussed in \cite{BW80} and \cite{BLS1975}.

Let $D$, $P$ and $Q$ be integers with $P > 0$ and $D = P^2 - 4Q \neq 0$.
Define $U_0 = 0$, $U_1 = 1$, $V_0 = 2$ and $V_1 = P$.
The Lucas sequences $U_k$ and $V_k$ with parameters $P$ and $Q$ are defined recursively for $k \ge 2$ by
\[
U_k = PU_{k-1} - QU_{k-2}  \quad \text{and} \quad V_k = PV_{k-1} - QV_{k-2} \thinspace .
\]
For $k \ge 0$ we also have
\[
U_k = (\alpha^k - \beta^k)/(\alpha - \beta) \quad \text{and} \quad
V_k = \alpha^k + \beta^k \thinspace ,
\]
where $\alpha$ and $\beta$ are the distinct roots of $x^2 - Px + Q = 0$.
Note that $\alpha \beta = Q$ and $\alpha + \beta = P$.

When $n$ is an odd positive integer, write $\delta(n) = n - (D/n)$ where $(D/n)$ is the Jacobi symbol.
It is known \cite[pp. 1391-1392]{BW80}, \cite[Theorem 8]{BLS1975} that if $n$ is prime and $(n, Q) = 1$, then
\begin{align}
U_{\delta(n)} & \equiv 0 \qmod{n}, \label{eq1} \\
V_{\delta(n)} & \equiv 2Q^{(1-(D/n))/2} \qmod{n}, \quad\mathrm{provided}\quad (n,D)=1, \label{eq2} \\
U_n & \equiv (D/n) \qmod{n}, \label{eq3} \\
V_n & \equiv V_1 = P \qmod{n}. \label{eq4}
\end{align}
If $(n,2PQD)=1$, any two of these congruences imply the other two.

Lucas pseudoprimes were defined in \cite{BW80}.
These are analogues of Fermat pseudoprimes in which $a^{n-1}-1$ is replaced by a Lucas sequence.

For reasons discussed in that paper, to use Lucas sequences for primality testing, we choose an algorithm for picking $D$, $P$, and $Q$ based on $n$, and we require that the Jacobi symbol $(D/n) = -1$.
If $n$ is prime, $(n, D) = (n, Q) = 1$, and $(D/n) = -1$, then $\delta(n) = n + 1$ and congruences \eqref{eq1} and \eqref{eq2} become \eqref{E:UCongruence} and \eqref{E:VCongruence}.
These two congruences are key parts of the primality test that we propose below.

We'll discuss \eqref{E:VCongruence} in detail later.
The other congruences, \eqref{eq3} and \eqref{eq4}, also hold if $n$ is prime, but these congruences are not very useful in primality testing \cite[Section 6]{BW80}:
most composite $n$ that satisfy congruence \eqref{eq3} have small prime factors;
many composite $n$ that satisfy \eqref{eq4} are psp(2).

If $n$ satisfies \eqref{E:UCongruence}, then $n$ is called a {\em Lucas probable prime} with parameters $P$ and $Q$, written lprp($P$, $Q$).
If $n$ satisfies \eqref{E:UCongruence} and we know it is composite, then we call $n$ a {\em Lucas pseudoprime}, written lpsp($P$, $Q$).
If $n$ fails \eqref{E:UCongruence}, then $n$ is composite.

For convenience, we also introduce the following

\medskip
\textbf{Definition} \cite[p.\ 266]{BressoudWagon}. If $n$ satisfies \eqref{E:VCongruence}, we call $n$ a \emph{Lucas-V probable prime} (vprp).
If $n$ is composite and satisfies \eqref{E:VCongruence} with parameters $P$ and $Q$, we call $n$ a \emph{Lucas-V pseudoprime} (vpsp$(P, Q)$).

\medskip
What we call vpsp's are sometimes called \emph{Dickson pseudoprimes of the second kind} \cite{Rot2003}.

The authors of \cite{BW80} proved that there are infinitely many Lucas pseudoprimes,
but that they are rare compared to the primes.

The precise sequence of numbers that turn out to be Lucas pseudoprimes depends on the algorithm for choosing $D$, $P$, and $Q$.
One algorithm, first proposed by John Selfridge in \cite{PSW1980} and mentioned in \cite{BW80}, and which seems to be widely used in primality testing, is:

\textbf{Method A}: Let $D$ be the first element of the sequence 5, $-7$, 9,
$-11$, 13, $-15$, $\ldots$ for which $(D/n)=-1$.
Let $P=1$ and $Q=(1-D)/4$.

This algorithm never sets $Q = 1$, but if $D = 5$, it sets $Q = -1$.
(Method A sets $Q = -1$ fairly often, namely, when $n \equiv \pm3 \qmod{10}$.)

We remarked in \cite{BW80} that more composite $n$ satisfying any
of \eqref{eq1}--\eqref{eq4} had $Q \equiv \pm1$ than $Q \not \equiv \pm1 \qmod{n}$.
This observation led the authors to define the following preferred method to select
parameters, which forces $Q \not \equiv \pm1 \qmod{n}$:

\textbf{Method A*}: Choose $D$, $P$, and $Q$ as in Method A above.
If $Q = -1$, change both $P$ and $Q$ to 5.

Method A* leaves $D = P^2 - 4Q$ unchanged from Method A.

It turns out that the Lucas pseudoprimes generated by Methods A and A* are the same.
The same is true for strong Lucas pseudoprimes (see Section \ref{S:IntroStrongLucas}).
We prove this in Appendix \ref{S:Astar}.

If $D$, $P$, and $Q$ are chosen with Method A*, the first ten lpsp are: 323, 377, 1159, 1829, 3827, 5459, 5777, 9071, 9179, and 10877.

Calculations performed for this paper show that when Method A* is used, there are 2\,402\,549 lpsp less than $10^{15}$.

\subsection{Strong Lucas probable primes and pseudoprimes} \label{S:IntroStrongLucas}

\cite{BW80} defines {\em strong} Lucas probable primes and {\em strong} Lucas pseudoprimes.
If $n$ is odd, then we can write $n + 1 = d \cdot 2^s$ where $d$ is odd.
If $n$ is prime and $(D/n) = -1$, then we will have either
\begin{align}
U_d             & \equiv 0 \qmod{n}, \quad\text{or} \label{E:StrongLrp1} \\
V_{d \cdot 2^r} & \equiv 0 \qmod{n}, \quad\text{for some } r \text{ with } 0 \le r < s. \label{E:StrongLrp2}
\end{align}

If $(D/n) = -1$ and $n$ satisfies \eqref{E:StrongLrp1} or \eqref{E:StrongLrp2}, then $n$ is called a {\em strong Lucas probable prime} with parameters $P$ and $Q$, written slprp($P$, $Q$).
If $n$ is an slprp($P$, $Q$), then $n$ is also an lprp($P$, $Q$), that is, $U_{n+1} = U_{d \cdot 2^s} \equiv 0 \qmod{n}$.

If $(D/n) = -1$, $n$ satisfies \eqref{E:StrongLrp1} or \eqref{E:StrongLrp2} and is {\em composite}, then $n$ is called a {\em strong Lucas pseudoprime}, written slpsp($P$, $Q$).

If $D$, $P$, and $Q$ are chosen with method A*, the first ten slpsp are: 5459, 5777, 10877, 16109, 18971, 22499, 24569, 25199, 40309, and 58519.

The slpsp($P, Q$) are scarcer than lpsp($P, Q$).
For example, of the 2\,402\,549 lpsp less than $10^{15}$, only 474971 are slpsp.

Because strong lpsp are rarer than lpsp, a sensible primality test will
use the strong version of the Lucas test, Congruences \eqref{E:StrongLrp1}
and \eqref{E:StrongLrp2}, instead of \eqref{E:UCongruence}.

The following equations show how to use the binary representation of $n+1$ to efficiently compute the values on the left sides of Congruences \eqref{E:StrongLrp1} and \eqref{E:StrongLrp2}.
We can also compute $U_{n+1}$, and, at almost no added computational cost, $V_{n+1}$ and $Q^{n+1}$.
\begin{align}
U_{2k} & = U_k V_k \label{eq5} \\
V_{2k} & = V_k^2 - 2Q^k \label{eq6} \\
Q^{2k} & = (Q^k)^2 \label{eq7} \\
U_{k+1} & = (PU_k+V_k)/2 \label{eq8} \\
V_{k+1} & = (DU_k+PV_k)/2 \label{eq9} \\
Q^{k+1} & = Q\cdot Q^k \label{eq10}
\end{align}
Equations \eqref{eq5} and \eqref{eq6} are Equations 4.2.6 and 4.2.7
in Williams \cite{Wil98} while \eqref{eq8} and \eqref{eq9} are 4.2.21
in that book.
Equations \eqref{eq5}--\eqref{eq7} are used to double the subscript
and exponent;
Equations \eqref{eq8}--\eqref{eq10} are used to increment the subscript and exponent.
These equations are also given in \cite[p. 628]{BLS1975}.

In Equations \eqref{eq8} and \eqref{eq9}, if the numerator is odd, we increment it by $n$ to make it be even.
This is legitimate because $n$ is odd, and we care only about the result modulo $n$.



For example, suppose $n = 323 = 17 \cdot 19$.
The Jacobi symbol $(5/n) = -1$, so Method A* sets $D = P = Q = 5$.
Writing $n + 1 = 324 = d \cdot 2^s$ where $d$ is odd, gives $d = 81$ and $s = 2$.
In binary, $d = 1010001$ and $n + 1 = 101000100$.
$n + 1$ has 9 bits, so there are (at most) 9 steps.
We process the bits from the left.
At each step, we double the subscript and exponent and compute $U_{2k}$, $V_{2k}$, and $Q^{2k} \qmod{n}$.
If the corresponding bit is odd, we also compute $U_{2k+1}$, $V_{2k+1}$, and $Q^{2k+1} \qmod{n}$.

All results are shown $\qmod{n}$. Step 1 merely initializes the sequences.
\begin{itemize}  
 \item bit 1 is 1: $(U_0, V_0, Q^0) = (0, 2, 1)$, $(U_1, V_1, Q^1) = (1, 5, 5)$;
 \item bit 2 is 0: $(U_2, V_2, Q^2) = (5, 15, 25)$;
 \item bit 3 is 1: $(U_4, V_4, Q^4) = (75, 175, 302)$, $(U_5, V_5, Q^5) = (275, 302, 218)$;
 \item bit 4 is 0: $(U_{10}, V_{10}, Q^{10}) = (39, 5, 43)$;
 \item bit 5 is 0: $(U_{20}, V_{20}, Q^{20}) = (195, 262, 234)$;
 \item bit 6 is 0: $(U_{40}, V_{40}, Q^{40}) = (56, 23, 169)$;
 \item bit 7 is 1: $(U_{80}, V_{80}, Q^{80}) = (319, 191, 137)$, $(U_{81}, V_{81}, Q^{81}) = (247, 306, 39)$;
 \item bit 8 is 0: $(U_{162}, V_{162}, Q^{162}) = (0, 211, 229)$;
 \item bit 9 is 0: $(U_{324}, V_{324}, Q^{324}) = (0, 135, 115)$;
\end{itemize}

From the last step, we see that $n$ is a Lucas probable prime because $U_{n+1} \equiv 0 \qmod{n}$.

$d = 81$, but none of $U_d$, $V_d$, or $V_{2d} = V_{162}$ is $0 \qmod{n}$.
Therefore, this $n$ fails the strong Lucas test, and we have proved that $n$ is composite upon reaching subscript $162 = (n+1)/2$.

Because of Congruences (\ref{E:StrongLrp1}) and (\ref{E:StrongLrp2}), the strong test always terminates at or before subscript $d \cdot 2^{s-1} = (n + 1)/2$
either with the conclusion that $n$ is composite, or that $n$ is a strong Lucas probable prime.

However, if $n$ \emph{is} a strong lprp, then we we recommend continuing the calculation a few more steps until we reach subscript $n + 1$, in order to obtain $V_{n+1} \qmod{n}$.
This enables us to check Congruence (\ref{E:VCongruence}), which, as we will see in Section \ref{S:data15}, is even more effective than (\ref{E:UCongruence}) at distinguishing primes from composites.

\section{The original Baillie-PSW primality test} \label{S:OriginalBPSW}

In \cite{BW80}, the first and third authors show that we get a very effective test for primality by combining Fermat and Lucas probable prime tests.

This combined test works so well because, in some sense, psp's and lpsp's tend to be different kinds of numbers.
For example, the numbers that are psp(2) and those that are lpsp from Method A* tend to fall into residue classes $+1$ and $-1$, respectively, for small moduli \cite[pp. 1404-1405]{BW80}.
A similar phenomenon is observed for psp($a$) for other $a$, and for lpsp's generated by several other methods for choosing $D$, $P$, and $Q$.

The probable prime test we proposed in \cite{BW80} has these steps:
\begin{enumerate}
\item If $n$ is not a strong base-2 probable prime, then $n$ is composite, so stop.
\item Choose Lucas parameters with Method A*.
(If you encounter a $D$ for which $(D/n)=0$: if either $|D| < n$, or if $|D| \ge n$ but $n$ does not divide $|D|$, then $n$ is composite, so stop.)
\item If $n$ is not a strong Lucas probable prime with the chosen parameters, then $n$ is composite.
Otherwise, declare $n$ to be (probably) prime.
\end{enumerate}

If $n$ is composite, the test almost always stops in the first step so the other steps are not needed.
The test almost never stops in the second step.
If $n$ is prime, then all three steps are needed.

The authors of \cite{PSW1980} and \cite{BW80} observed that, up to $25 \cdot 10^9$, there was no overlap between the psp(2) and the lpsp from Method A*.
Using more recent data from Feitsma \cite{Feitsma}, we find that none of the 118\,968\,378 psp(2) up to $2^{64} \approx 1.8 \cdot 10^{19}$ is an lpsp when Method A* is used.
Therefore, this test correctly distinguishes primes from composites up to at least $2^{64}$.
Further, no one has reported a larger composite $n$ that is both psp(2) and lpsp($P$, $Q$) using method A*.

Richard Pinch \cite{Pinch}, \cite{Pinch18} has computed a list of all 20\,138\,200 Carmichael numbers up to $10^{21}$.
He kindly provided his list to the first author, for which we thank him.
None of these Carmichael numbers is an lpsp when Method A* is used.

The reader might notice that the above test does not first check $n$ for divisibility by small primes.
This check is omitted because it is not necessary (although step 2 does sometimes find small factors).
However, \emph{for the sake of efficiency}, a practical primality test should first check to see whether $n$ is divisible by small primes before proceeding to step 1.

Some version of this test is used as a fast algorithm for finding large probable primes
in mathematical software packages like \textit{FLINT}, \textit{Maple}, \textit{Mathematica}, \textit{Pari/GP}, \textit{SageMath}, and by
programs for choosing large primes for public-key ciphers like RSA.

Several programming languages, like \textit{GNU GMP}, \textit{Java}, and \textit{Perl} also provide functions for doing Fermat and Lucas tests.

While some cryptographic libraries use a combined Fermat/Lucas test, some do not.
Albrecht, et al, were able to find composite numbers which some of the latter libraries declared were prime \cite{Albrecht}.

A reward of \$620 was offered for an example of a composite $n$
declared prime by this test.
No one has claimed the reward after 40 years;
many have tried to collect it.
It has been tested on billions of large odd integers $n$
and has never been reported to have failed.

\section{The data to \texorpdfstring{$10^{15}$}{1E15}}\label{S:data15}

Recall that lpsp and vpsp are composite $n$ that satisfy \eqref{E:UCongruence} and \eqref{E:VCongruence}, respectively.

We computed the lpsp and vpsp up to $10^{15}$, using Method A* to choose the Lucas parameters.

This calculation took about 750000 core-hours on the Rice cluster at Purdue University,
plus about 10000 core-hours on computers at the University of Lethbridge.

The counts are shown in Table \ref{Ta:AllCountsTo10k}.
What is striking is that, while there are about 2 million each of psp(2) and lpsp, there are only \emph{five} vpsp.

These five numbers are shown in Table \ref{Ta:VSolnCounts}.
For $n = 14\,760\,229\,232\,131$, Method A* set $P = 1$, $Q = 2$.
For the other four $n$, Method A* set $P = Q = 5$.

\begin{table}[ht]
\caption{Number of psp(2), spsp(2), lpsp, slpsp, and vpsp with $n < 10^k$ using Method A*.}
\label{Ta:AllCountsTo10k}
\begin{tabular}{|r|r|r|r|r|r|} \hline
$k$ & psp(2) & spsp(2) & lpsp & slpsp & vpsp \\ \hline
 2 &       0 &      0 &       0 &      0 & 0 \\
 3 &       3 &      0 &       2 &      0 & 1 \\
 4 &      22 &      5 &       9 &      2 & 1 \\
 5 &      78 &     16 &      57 &     12 & 1 \\ \hline
 6 &     245 &     46 &     219 &     58 & 1 \\
 7 &     750 &    162 &     659 &    178 & 1 \\
 8 &    2057 &    488 &    1911 &    505 & 1 \\
 9 &    5597 &   1282 &    5485 &   1415 & 1 \\
10 &   14884 &   3291 &   15352 &   3622 & 1 \\ \hline
11 &   38975 &   8607 &   42505 &   9714 & 1 \\
12 &  101629 &  22407 &  116928 &  25542 & 3 \\
13 &  264239 &  58892 &  319687 &  67045 & 3 \\
14 &  687007 & 156251 &  875270 & 178118 & 4 \\
15 & 1801533 & 419489 & 2402549 & 474971 & 5 \\ \hline
\end{tabular}
\end{table}

\begin{table}[ht]
\caption{vpsp $<10^{15}$ using Method A*.}
\label{Ta:VSolnCounts}
\begin{tabular}{|r|l|l|l|} \hline
$n$ & $n$ factored & $n-1$ factored & $n+1$ factored \\ \hline
913 & $11 \cdot 83$ & $2^4 \cdot 3 \cdot 19$ & $2 \cdot 457$ \\  \hline
150\,267\,335\,403 & $3 \cdot 47 \cdot 89 \cdot 563 \cdot 21269$ & $2 \cdot 157 \cdot 478\,558\,393$ & $2^2 \cdot 1609 \cdot 23\,347\,939$ \\  \hline
430\,558\,874\,533 & $75913 \cdot 5\,671\,741$ & $2^2 \cdot 3^2 \cdot 11\,959\,968\,737$ & $2 \cdot 197947 \cdot 1\,087\,561$ \\  \hline
14\,760\,229\,232\,131 & $2467 \cdot 5\,983\,068\,193$ &
\parbox{1.1in}{$2\cdot3\cdot5\cdot107\cdot53569\cdot\mathstrut$ $\quad\cdot85837$\vspace{1mm}}
& $2^2 \cdot 3\,690\,057\,308\,033$ \\  \hline
936\,916\,995\,253\,453 & $2027 \cdot 21521\cdot 21\,477\,559$ &
\parbox{1.1in}{$2^2\cdot3\cdot37\cdot41\cdot1109\cdot$ $\quad\cdot 46\,409\,057$\vspace{1mm}} &
\parbox{1.1in}{$2\cdot389\cdot15313\cdot\mathstrut$ $\quad\cdot 78\,643\,211$\vspace{1mm}}  \\ \hline
\end{tabular}
\end{table}

In Section \ref{S:Qgt1}, we give heuristic arguments as to why vpsp are so rare, especially when $Q \not \equiv \pm1 \qmod{n}$.

When $P$ and $Q$ are chosen by Method A*, we found that (i) none of the five vpsp$(P, Q)$ is an lpsp$(P, Q)$, (ii) none of the 118\,968\,378 psp(2) less than $2^{64}$ is either an lpsp or a vpsp, and (iii) none of the 20\,138\,200 Carmichael numbers below $10^{21}$ is either an lpsp or a vpsp.
(We do not know if there is an $n > 10^{15}$ which is both lpsp$(P, Q)$ and vpsp$(P, Q)$.)

The enhanced primality test we propose in Section \ref{S:newTest} is based on the rarity of vpsp, and on this absence of overlap between any two of spsp(2), slpsp, and vpsp.

Dana Jacobsen's website \cite{Jacobsen} displays counts of psp(2), spsp(2), lpsp, and slpsp less than $10^{15}$, where $P$ and $Q$ are selected by method A (or A*), as well as other types of pseudoprimes.
The lists of Lucas pseudoprimes can be downloaded from that site, or from the ancillary files that accompany version 1 of this preprint, \cite{BFWPreprint}.

\section{Is there anything special about Method A*?} \label{S:IsAStarSpecial}
The reader may wonder whether the rarity of vpsp compared to lpsp is an artifact of using Method A* to choose $D$, $P$, and $Q$.

The answer appears to be ``no'', especially if we require that $Q \not \equiv \pm1 \qmod{n}$.
We compared several methods for choosing $D$, $P$, and $Q$; see Table \ref{Ta:OtherMethods}.

For example, \cite{BW80} describes Methods B and B*:

\medskip
\textbf{Method B}: Let $D$ be the first element of the sequence 5, 9, 13, 17, $\ldots$ for which $(D/n) = -1$.
Let $P$ be the smallest odd number exceeding $\sqrt{D}$, and $Q = (P^2 - D)/4$.

\textbf{Method B*}: Choose $D$, $P$, and $Q$ as in Method B.
If $Q = 1$, replace $Q$ by $P+Q+1$ and replace $P$ by $P+2$ (this preserves the value of $D$).
\medskip

For this paper, we also tested:

\textbf{Method C}: Same as Method A, except we start testing $D$'s at $D = 41$ instead of at $D = 5$.
This insures that Method C almost always produces a $(P, Q)$ pair different from the pair produced by Method A.

\textbf{Method D}: Fix Q = 2.
Try $P = 4, 5, 6, 7$, ... until $(D/n) = -1$.

\textbf{Method R1}: Choose $P$ and $Q$ \emph{at random} from a uniform distribution with $1 \le P, Q \le n - 1$, until $(D/n) = -1$.
We used the \emph{random( )} function in version 2.11.4 of \emph{PARI/GP}, initialized with \emph{PARI}'s default seed of 1.

\textbf{Method R2}: Same as Method R1, but initialized with the (randomly-selected) seed 737984.
\medskip

We compared these eight methods for odd, composite $n < 10^{10}$.

Methods A*, B*, C, and D can never set $Q \equiv \pm1 \qmod{n}$.

Method B yielded 5940 vpsp.
Only one of these, $n = 64469$, occurred with $Q \not \equiv \pm1 \qmod{n}$: This $n$ is vpsp$(5, 3)$, but is not lpsp$(5, 3)$.

Method B* yielded two vpsp: $n = 913$ ($P = Q = 5$) and $n = 64469$ ($P = 5$, $Q = 3$).

No vpsp from Method R1 or R2 had $Q \equiv \pm1 \qmod{n}$.
This is not surprising: $Q$ was a random integer between 1 and $n-1$, and $Q \equiv \pm1 \qmod{n}$ occurred for only eight $n$'s with R1 and twelve with R2.

\begin{table}[ht]
\begin{tabular}{|l|r|r|r|r|} \hline
 & & & & simultaneously \\ Method & lpsp & vpsp, $Q \equiv \pm1$ & vpsp, $Q \not \equiv \pm1$ & lpsp and vpsp \\ \hline
 A  &   15352 &           914 &                    0 &                     757 \\
 A* &   15352 &           --  &                    1 &                       0 \\
 B  &   15019 &          5939 &                    1 &                    4374 \\
 B* &   12879 &           --  &                    2 &                       0 \\
 C  &   13766 &           --  &                    4 &                       0 \\
 D  &   15957 &           --  &                    6 &                       0 \\  
 R1 &   17065 &             0 &                    3 &                       0 \\
 R2 &   16863 &             0 &                    4 &                       0 \\ \hline
\end{tabular}
\caption{Number of lpsp and vpsp to $10^{10}$ using various methods for choosing $P$ and $Q$.}
\label{Ta:OtherMethods}
\end{table}

All of the methods tested in Table \ref{Ta:OtherMethods} yielded fewer vpsp than lpsp.
Moreover, if $Q \not \equiv \pm1 \qmod{n}$, none of the lpsp$(P, Q)$ was also vpsp$(P, Q)$.
This Table supports the importance of choosing $Q \not \equiv \pm1 \qmod{n}$.

\section{The enhanced BPSW primality test} \label{S:newTest}

The enhanced primality test we propose here is based on the one described in Section \ref{S:OriginalBPSW}.
The most important strengthening is that we now include Congruence \eqref{E:VCongruence} to check whether $n$ is a vprp.
This has very little additional computational cost beyond the Lucas test in step 3.

The strong Lucas probable prime test, Congruences
\eqref{E:StrongLrp1} and \eqref{E:StrongLrp2}, allows us to stop
the calculation one or more steps before reaching
$U_{n+1}$, $V_{n+1}$, and $Q^{n+1} \qmod{n}$.
Here, we assume that we continue the calculation for a few additional steps in order to obtain $Q^{(n+1)/2}$ and $V_{n+1} \qmod{n}$.

Here is our proposed enhanced primality test for odd, positive integer $n$:
\begin{enumerate}
\item If $n$ is not a strong probable prime to base $2$, then $n$ is composite; stop.
\item Choose Lucas parameters with Method A*.
If you encounter a $D$ for which $(D/n)=0$: if either $|D| < n$, or if $|D| \ge n$ but $n$ does not divide $|D|$, then $n$ is composite; stop.
\item If $n$ is not an slprp($P$, $Q$), then $n$ is composite; stop.
\item If $n$ is not a vprp($P$, $Q$), then $n$ is composite; stop.
\item If $n$ does not satisfy $Q^{(n+1)/2} \equiv Q \cdot (Q/n) \qmod{n}$, then $n$ is composite; stop.
Otherwise, declare $n$ to be probably prime.
\end{enumerate}

Recall that no composite number is known that passes steps 1 through 3.
This test is more powerful than the original BPSW test because so few composite $n$ satisfy step 4.
A composite $n$ that passes this test would have to be, simultaneously, spsp(2), slpsp$(P, Q)$, \emph{and} vpsp$(P, Q)$.
Consequently, we expect that a composite $n$ would be even less likely to pass this test than to pass the original BPSW test.

An odd, composite $n$ that is both an lpsp and a vpsp is a \emph{Frobenius}
pseudoprime \cite[p. 145]{CP}, \cite{GranthamHC}, \cite{GranthamFP}.
These are rare \cite{Jacobsen}, in part because, as we've seen above, the vpsp are rare.
Odd, composite $n$ that pass this enhanced test should be even rarer.

Step 5 is a primality check based on Euler's criterion, Congruence \eqref{E:EulerCriterion}.
This is a relatively minor enhancement.
However, since we essentially already have the power of $Q$ necessary for the test, we may as well use it.
Once we have calculated $Q^{(n+1)/2} = Q \cdot Q^{(n-1)/2} \qmod{n}$, we can compute $(Q/n)$, then apply Euler's criterion to check whether
\begin{equation*}
Q^{(n+1)/2}\equiv Q \cdot (Q/n) \qmod{n} \thinspace .
\end{equation*}
If this congruence fails, then $n$ is composite.

\medskip

\textbf{Suggestions for implementing this primality test.}

1. For efficiency, before step 1, one should first check $n$ for divisibility by small primes.

2. We recommend doing a (strong) Fermat test to base 2 instead of to some other base.
As far as anyone knows, there is nothing inherently better about using base 2.
However, because we know all psp(2) up to $2^{64}$, we know that no psp(2) below that limit is an lpsp.
We do not know whether this is true for other bases.

3. In step 2: If $n$ happens to be a perfect square, then $(D/n)$
will never be $-1$. So, after encountering, say, 20 $D$'s with $(D/n) = 1$, one
should check whether $n$ is a perfect square; if so, it is composite.
This can be done quickly using Newton's method; see \cite[p. 1401]{BW80}.

4. It is easy to show that, if $n$ is sprp($a$), then $n$ is also sprp($\pm(a^k)$) for $k \ge 1$.
Therefore, if Method A* chooses a $Q$ such that $|Q|$ is a power of 2, then, because $n$ is
known from step 1 to be sprp(2), the test in step 5 will not strengthen the test.
For $n < 10^9$, this happens about 28 percent of the time ($D = -7$, $Q = 2$; $D = -15$, $Q = 4$; $D = 17$, $Q = -4$, etc).

5. To compute all three of $U_{n+1}$, $V_{n+1}$, and $Q^{n+1}$ $\qmod{n}$ takes roughly three times as many multi-precision operations as it takes to compute $2^{n-1} \qmod{n}$.
Therefore, this enhanced BPSW test takes about as long as doing Fermat tests to four different bases.
However, as noted in \cite[p. 1020]{PSW1980}, if $n$ is psp to base $a$, then $n$ is more likely than the average number of that size to also be psp to some other base $b$.
In other words, there are diminishing returns in doing repeated Fermat tests.
Therefore, it makes more sense to do one Fermat test followed by the Lucas tests in steps 3 and 4, than to perform Fermat tests to four (or more) different bases.

\medskip

\textbf{Reward for a counterexample or for a proof that there are none.}

A counterexample to this enhanced test would be a positive, odd composite $n$ which this test declares is probably prime.
The first and third authors each offer U.S. \$1000 for either the first counterexample to this enhanced test, or the first proof, published in a peer-reviewed journal, that there are none.
A claim that $n$ is a counterexample must be accompanied by a $(P, Q)$ pair that came from Method A*, for which the test claims $n$ is probably prime.
The claim must also be accompanied by a proof that $n$ is composite: either a (not necessarily prime) factor of $n$ that is larger than 1 and less than $n$,
or a base $a$ with $2 < a < n-1$ for which $n$ is not a base-$a$ strong probable prime,
or a $(P, Q)$ pair not from Method A* such that $D = P^2 - 4Q$ has Jacobi symbol $(D/n) = -1$, but for which $n$ is not slpsp$(P, Q)$, or is not vpsp$(P, Q)$.

\section{Some Heuristics for lpsps and vpsps} \label{S:Qgt1}

We note that both conditions of \eqref{eq1} and \eqref{eq2} are both congruences modulo $n$, and thus
if $p$ is any prime which divides $n$, then we obtain implied congruences modulo $p$.
That is if
\[ U_{\delta(n)} \equiv 0 \qmod{n}, \qquad \textrm{then}\qquad U_{\delta(n)} \equiv 0 \qmod{p} \]
and if 
\[ V_{\delta(n)}  \equiv  2Q^{(1-(D/n))/2} \qmod{n}, \qquad\textrm{then}\qquad V_{\delta(n)}  \equiv  2Q^{(1-(D/n))/2} \qmod{p}. \]
Moreover, if $n$ were square free, we have that the conditions
modulo $p$ for all $p$ dividing $n$ would give sufficient conditions
for the same congruences modulo $n$.

Now, suppose we write $n=py$, we can assess the probability that for example
\[ V_{\delta(n)}  \equiv  2Q^{(1-(D/n))/2} \qmod{p} \]
by considering the probability that, as we vary $y$ among $y$ with $(D/py) = (D/n)$,
\[ V_{\delta(py)}  \equiv  2Q^{(1-(D/py))/2} \qmod{p}. \]
We note that this quantity only depends on $y$ modulo the period of the sequence $V_k \qmod{p}$, so the probability is well defined.
If these probabilities were independent for distinct prime
factors of $n$, then for square free $n$ we could determine the
probability that $n$ is a pseudoprime from local contributions.
These probabilities are most likely not independent.

Assume D is square free.

Before looking at these probabilities in the various cases we first recall a few facts.
Fix the prime $p$.
Let ${\mathbb F}_p$ and ${\mathbb F}_{p^2}$ denote the fields
of order $p$ and $p^2$, respectively.
In the notation of Section \ref{S:IntroLucas},
$$\alpha=\frac{P+\sqrt{D}}{2}\qquad \beta=\frac{P-\sqrt{D}}{2}.$$
We think of these quantities as elements of ${\mathbb Q}(\sqrt{D})$,
or ${\mathbb F}_{p^2}$ when $(D/p)=-1$, or ${\mathbb F}_p$ when $(D/p)=1$.

When $(D/p)=-1$, for $x$, $y\in{\mathbb F}_p$ we have
$(x+y\sqrt{D})^p=x-y\sqrt{D}$ in ${\mathbb F}_{p^2}$.
In particular, $\alpha^p=\beta$ and $\beta^p=\alpha$.
Whereas when $(D/p)=1$, for $x$, $y\in{\mathbb F}_p$ we have
$(x+y\sqrt{D})^p=x+y\sqrt{D}$ and $x^{p-1}=1$ in ${\mathbb F}_p$.
In particular, $\alpha^{p-1}=1$ and $\beta^{p-1}=1$.

These facts allow us to derive the formulas \eqref{eq1} and \eqref{eq2}.
Indeed if $(D/p)=-1$, then in ${\mathbb F}_{p^2}$ we have
$$U_{p+1}=\frac{\alpha^{p+1} - \beta^{p+1}}{\alpha - \beta} =
\frac{\alpha\beta-\beta\alpha}{\alpha - \beta} = 0$$
and 
$$V_{p+1} = \alpha^{p+1} + \beta^{p+1} =2\alpha\beta=2Q.$$
Since both sides are in ${\mathbb F}_p$, we may think of
these congruences as congruences modulo $p$.
These give exactly \eqref{eq1} and \eqref{eq2} with $(D/p)=-1$.

Similarly if $(D/p)=1$, then in ${\mathbb F}_p$ we have
$$U_{p-1}=\frac{\alpha^{p-1} - \beta^{p-1}}{\alpha - \beta} =
\frac{1-1}{\alpha - \beta} = 0$$
and 
$$V_{p-1} = \alpha^{p-1} + \beta^{p-1} =1+1=2.$$

We are most interested in the case $(D/n)=-1$ and composite $n$
that satisfy $U_{n+1}\equiv 0\qmod{n}$ and $V_{n+1} \equiv 2Q\qmod{n}$.
Consequently in the following we shall focus on the case of $(D/n)=-1$.
We now investigate the probabilities mentioned above.

\subsection{Case 1: \texorpdfstring{$(D/p)=-1$}{(D/p)=-1}}\label{case1}

Since $(D/p)(D/y)=(D/py)=(D/n)=-1$, we must have $(D/y)=1$.
We also have
$$U_{py+1}=\frac{\alpha^{py+1} - \beta^{py+1}}{\alpha - \beta} =
\alpha\beta \frac{\beta^{y-1}-\alpha^{y-1}}{\alpha - \beta} = QU_{y-1}$$
and
$$V_{py+1}=\alpha^{py+1} + \beta^{py+1} =
Q(\beta^{y-1}+\alpha^{y-1})=QV_{y-1}.$$
Thus, for Congruences \eqref{eq1} and \eqref{eq2} we are interested
respectively in the probability that $U_{y-1}\equiv0\qmod{p}$ and
the probability that $V_{y-1}\equiv2\qmod{p}$.

The sequences $U_y$ and $V_y$ as functions of $y$ modulo $p$
are periodic with periods less than $p^2$.
The condition $(D/y)=1$ has period $D$ (or $4D$) as we are implicitly interested in representatives for $y$ that are odd. By the CRT this modulo $D$ condition is irrelevant to the conditional probability unless $D\mid p^2-1$.

\begin{lemma}\label{lemma4}
The periodicity of the appearance of 0 for the sequence $U_{y-1}$
is exactly the order of the image of $\alpha$ in the group
${\mathbb F}_{p^2}^{\times}/{\mathbb F}_{p}^{\times}$.
In particular, it divides $p+1$.

Moreover, the order is 2 when $P=0 \qmod{p}$.
\end{lemma}

\noindent
{\em Proof}.
We recall that the sequence $U_{y-1}$ gives the irrational part of the number $\alpha^{y-1}$, and thus
$U_{y-1} = 0 \qmod{p}$ if and only if $\alpha^{y-1} \in {\mathbb F}_{p}$.
Thus the collection of $y$ such that $U_{y-1}=0$ is exactly the
kernel of the map
${\mathbb F}_{p^2}^{\times}\rightarrow
{\mathbb F}_{p^2}^{\times}/{\mathbb F}_{p}^{\times}$.
Hence, we are really studying the image of the map
$\langle \alpha^n \rangle \rightarrow
{\mathbb F}_{p^2}^{\times}/{\mathbb F}_{p}^{\times}$.
Since the image is a subgroup, its order divides $p+1$.
This completes the proof of the first claim.

We note that if $P=0 \qmod{p}$,
then $\alpha = \sqrt{D}/2$ and $\alpha^2=-Q\in {\mathbb F}_{p}$.

Since we know $U_{(p+2)-1}=0$, we obtain
\begin{proposition}
There exists a divisor $k$ of $p+1$ such that $U_{py+1}\equiv0\qmod{p}$
if and only if $y\equiv p+2 \equiv 1\qmod{k}$.

In particular, the probability that $U_{py+1}\equiv0\qmod{p}$ is $\frac{1}{k} \ge \frac{1}{p+1}$.

If we assume additionally that $y$ is odd then subject to this condition the probability is at least $ \frac{2}{p+1}$.
\end{proposition}
We note that in the above we are not directly accounting for the possibility
that $D$ divides $k$, in which case 
the condition $(D/y)=1$ would tend to lead to a higher conditional probability as some of the $k$ options for $y$ must be discarded, but we will never discard the option $1$ mod $k$.
In particular the claim remains valid with these considerations.

Now we note that the order of $\alpha$ in ${\mathbb F}_{p^2}^{\times}/{\mathbb F}_{p}^{\times}$ is precisely the order
of $\frac{\alpha}{\beta} = \frac{\alpha^2}{Q}$ in ${\mathbb F}_{p^2}^{\times}$ and the condition $U_{y+1}\equiv 0\qmod{p}$ is 
equivalent to $ \left(\frac{\alpha^2}{Q}\right)^{y+1} = 1 \in {\mathbb F}_{p^2}^{\times}$.
In the case $n$, and hence $y$, are odd we will always have $y+1$ even, and writing $y+1=2z$ the condition on $z$ is
\[  1 = \left( \frac{\alpha^2}{Q}\right)^{2z}  . \]
As such we can see that this order is automatically at most $(p+1)/2$.

If we continue to write $y+1=2z$ and look at
\[  1 = \left( \frac{\alpha^2}{Q}\right)^{2z}  =  \left( \frac{\alpha^2}{-Q}\right)^{2z} . \]
we can see that the $2$-part of the order of $ (\alpha^2/Q)^2$ will tend to further bounded if $Q$ or $-Q$ is a square. 
Though this is guaranteed when $p\equiv 3\qmod{4}$ it is also guaranteed if $Q$ or $-Q$ is an integer perfect square. 

Additionally we notice that when $y$ is odd if we replace $\alpha,\beta$ by the conjugate pair $\alpha' = \sqrt{D}\alpha, \beta' = -\sqrt{D}\beta$ we effectively replace 
 $Q$ by $-DQ$ but the order of $(\alpha'/\beta')^2$ agrees with that of $(\alpha/\beta)^2$.
It follows that having any of 
\[ Q,\quad, -Q,\quad DQ, \quad -DQ \]
perfect integer squares will tend to increase the probability $U_{y+1}\equiv0\qmod{p}$ when $y$ is odd. This phenomenon can be observed empirically by
counting the proportion of $n$ which are lpsp for different options $P,Q$.
Those with $Q$ of the above form tend to appear more than on average.

The situation for \eqref{eq2} is somewhat more subtle. We have the following
\begin{lemma}\label{lemma8}
The period of the sequence $V_{y-1}$ divides the order of $\alpha$ in
${\mathbb F}_{p^2}^{\times}$, which is a divisor of $p^2-1$.

Moreover, the order of $\alpha$ in
${\mathbb F}_{p^2}^{\times}$ is divisible by the LCM of the order
of $\alpha$ in the group 
${\mathbb F}_{p^2}^{\times}/{\mathbb F}_{p}^{\times}$
and the order of $Q$ in ${\mathbb F}_{p}^{\times}$, and is at most twice this amount.
In particular, this period is at least as large as the period
of Lemma \ref{lemma4}.
\end{lemma}

\noindent
{\em Proof}.
The first claim about the period is clear given that $\alpha$ and $\beta$ have the same period.

For the second claim we note that the order of an element is divisible by the order of its image under any homomorphism.
We obtain the result by considering the map
${\mathbb F}_{p^2}^{\times} \rightarrow {\mathbb F}_{p^2}^{\times}/{\mathbb F}_{p}^{\times}$ as well as the norm map
\[ N_{{\mathbb F}_{p^2}^{\times}/{\mathbb F}_{p}^{\times}} : {\mathbb F}_{p^2}^{\times} \rightarrow {\mathbb F}_{p}^{\times} \]
for which we have
\[ N_{{\mathbb F}_{p^2}^{\times}/{\mathbb F}_{p}^{\times}}(\alpha) = \alpha^{p+1} = Q \]
Because the intersection of the kernels of these two maps is $\pm1$ we conclude that the exact order of $\alpha$ is either the LCM or the two quantities, or exactly twice this.
(The exact order is twice this if and only if there exists $z$ with $\alpha^z = -1 \qmod{p}$. This is guaranteed if $Q$ is not a square mod $p$.)
We note further that because we are taking the LCM of a number dividing $p+1$ and one dividing $p-1$ the LCM is almost exactly the product.

Note from the lemma above we may conclude that if $P=0$ the order of $\alpha$ is exactly $2$ times the order of $-Q$.
We also conclude that if $Q=\pm1$, the order of $\alpha$ divides $2(p+1)$.

It remains the case that we cannot expect that $2$ only appears once in each period. 
Additionally, in contrast to the previous case there is no guarantee that the appearance of $2$ in the period is actually itself periodic.

\begin{lemma}
Let $\ell$ denote the period of $V_y$ modulo $p$.
 We have that $V_{y+m(p+1)} = Q^mV_y$ and consequently 
 each $a\in \mathbb{F}_p^\times$ is repeated by $V_y$ equally often as $Q^ma$ and hence not more than $\ell/\ord(Q)$ times within one period of $V_y\qmod{p}$.
\end{lemma}

\noindent
{\em Proof}.

This follows from the observation that
\[ V_{y + p+1} = QV_{y} \]
from which we obtain a bijection between the occurrences of $x$ and $Qx$.
Hence, each value $ 0\not\equiv x\qmod{p}$ which occurs, does so just as often as $Qx\qmod{p}$ 
in one period.

Now we consider the map
\[ \Psi : \mathbb{Z} \rightarrow {\mathbb F}_p \times {\mathbb F}_{p}^{\times} \]
given by
\[ y \mapsto  (V_y, N_{{\mathbb F}_{p^2}^{\times}/{\mathbb F}_{p}^{\times}}(\alpha^y) )=  (V_y,Q^y). \]
And note that this map has a period which is either $\ord(\alpha)$ or $\ord(\alpha)/2$

\begin{lemma}
  The function $\Psi$ is exactly $2:1$ on its image, and hence each $a\in \mathbb{F}_p^\times$ is repeated by $V_y$ no more than $2\ord(Q)$ times in $0\leq y \leq {\ord(\alpha)}$.
\end{lemma}

We note that the image of $\Psi$ gives the trace and norm of $\alpha^y$, hence we can recover the minimal polynomial of $\alpha^y$ from $\Psi(y)$.
It follows that $\Psi(y_1) = \Psi(y_2)$ implies either $\alpha^{y_1} = \alpha^{y_2}$ or $\alpha^{y_1} = \beta^{y_2}$.
Because $\ell$ divides the order of $\alpha$ in the first case we obtain $y_1\equiv y_2 \qmod{\ell}$.
In the second case we obtain $y_1 \equiv py_2 \equiv -y_2 \qmod{\ell}$.

\begin{proposition}
The probability that $V_{py+1}\equiv 2Q \qmod{p}$ is less than the minimum of 
\[ \frac{1}{\ord(Q)} \qquad\textrm{and}\qquad \frac{2\ord(Q)}{\ord(\alpha)}. \]
We remark that 
\[
  \frac{2\ord(Q)}{\ord(\alpha)} =
  \begin{cases}
    \frac{1}{k} & \textnormal{if } \exists z,\alpha^z = -1 \qmod{p} \\
    \frac{2}{k} & \textnormal{ otherwise}
  \end{cases}
\]
where as before, $k$ is the order of $\alpha$ in ${\mathbb F}_{p^2}^{\times}/{\mathbb F}_{p}^{\times}$
\end{proposition}

We first note that this proposition gives an indication of why having $\ord(Q)$ large is beneficial.
Moreover, it indicates why one should expect the $V$ test to be better than the $U$ test, and at least as good even when $Q=-1$.

We next note that in the above we are not accounting for the possibility that $D \mid  \ell$, in which case 
the condition $(D/y)=1$ and $(D/p)=-1$ would imply the map $\Psi$ must be injective on relevant cases. 
The effect is that the conditional probabilities are still bounded by the above.

Finally we note that when we combine both the $U$ and $V$ conditions, and consider the conditional probability of the $V$ condition assuming the $U$ condition this amounts to restricting to the subsequence $y=(p+2)+kx$
for which that $U_{y-1}=0$, then as $(\alpha-\beta)U_z+2\beta^z=V_z$,
the condition $V_{y-1}=V_{(p+2)+xk-1}=2$ becomes
\[  \beta^{(p+1)+kx} = Q(\beta^{k})^x. \]
This condition is periodic in $x$, with period the exact order of $\beta^k$, as $k$ is the smallest power for which $\beta^k \in {\mathbb F}_{p}^{\times}$ the order
of $\beta^k$ divides $p-1$ and is, up to a multiple of $2$, the order of $Q$.
In particular, if $Q$ has a large order, the probability that the $V$ condition is satisfied remains low independently of the $U$ condition.

\subsection{Case 2: \texorpdfstring{$(D/p)=1$}{(D/p)=1}}\label{case2}

Since $(D/p)(D/y)=(D/py)=(D/n)=-1$, we must have $(D/y)=-1$.
As $\alpha^p=\alpha$ and $\beta^p=\beta$ we also have
$$U_{py+1} = \frac{\alpha^{py+1} - \beta^{py+1}}{\alpha-\beta} =
\frac{\alpha^{y+1} - \beta^{y+1}}{\alpha-\beta} = U_{y+1},$$
so we want to estimate the fraction of $y$ with $({D/y})=-1$ that have
$U_{y+1}\equiv0\qmod{p}$.
Likewise we have
$$V_{py+1} = \alpha^{py+1} + \beta^{py+1}=
 \alpha^{y+1} + \beta^{y+1} = V_{y+1},$$
so we want to estimate the fraction of $y$ with $({D/y})=-1$ that have
$V_{y+1}\equiv2Q\qmod{p}$.

We note that in this case $\alpha,\beta\in {\mathbb F}_p^\times$ are essentially independently chosen elements (determined by $P$ and $Q$). 

\begin{lemma}
The sequence $U_{y+1}$ is zero precisely when $((\alpha)\beta^{-1})^{y+1}=1$.
Hence the period of the vanishing of $U_{y+1}$ is precisely the
order of $(\alpha)\beta^{-1}$ as an element of ${\mathbb F}_{p}^{\times}$.
In particular, it divides $p-1$ and the LCM of the orders of
$\alpha$ and $\beta$.
Consequently, there is a divisor $k$ of $p-1$ such that $U_{y+1}=0$ if and only if $y=-1 \qmod{k}$.
\end{lemma}

\noindent
{\em Proof}.
The condition $U_{y+1}\equiv0\qmod{p}$ becomes
$\alpha^{y+1}=\beta^{y+1}$ in ${\mathbb F}_{p}^{\times}$, or
$(\alpha)\beta^{-1}=1$, which proves the lemma.

If $D\mid k$, then we should consider the impact of the condition $(D/y)=-1$. In contrast to the previous case 
it may not be possible to have $(D/y)=-1$ and  $y=-1 \qmod{k}$. If we assume $D\mid k$ so that $y=-1\qmod{D}$, then in the case
\begin{itemize}
\item $D=1\qmod{4}$ and $D>0$ then $(D/y)=(y/D)=(-1/D)=1$ hence the conditions are never simultaneously satisfiable.
\item $D=3\qmod{4}$ and $D>0$ then $(D/y)=(-1/y) (y/D)=(-1/y) (-1/D)=-(-1/y)$ hence the condition is satisfiable if $y=1\qmod{4}$. But we note that if $4\mid k$ this is not possible.
\item $D=3\qmod{4}$ and $D<0$ then $(D/y)=(-1/y) (y/-D)=(-1/y) (-1/-D)=(-1/y)$ hence  the condition is satisfiable if  $y=3 \qmod{4}$.
\item $D=1\qmod{4}$ and $D<0$ then $(D/y)=(y/-D)=(-1/-D)=-1$ hence the condition is always satisfiable.
\end{itemize}
The above suggests $D=1\qmod{4}$ and $D>0$ would be ideal.
We note about the above conditions, $D\mid p-1$ is only particularly likely for random $n$ when $D$ is small. Counts of lpsps are consistent with this expectation
that small positive $D$ have fewer lpsps.

We remark that as in the case $(D/p)=-1$ if $y+1=2z$  we are considering the condition
\[ 1 = \left(\frac{\alpha}{\beta}\right)^{y+1} = \left(\frac{\alpha^4}{Q^2}\right)^{z}. \]
In the case $p\equiv 3\qmod{4}$ the $2$-part of the order is already reduced to $1$.
However, if $p\equiv 1\qmod{4}$, we will have $-1,D,-D$ are all squares modulo $p$, hence if any of
\[ Q,\quad -Q, \quad DQ,\quad -DQ \]
are perfect integer squares this will reduce the maximum $2$-part of the order of $((\alpha)\beta^{-1})^2$ and hence increase the probability that $ U_{y+1}\equiv0\qmod{p}.$
In contrast if $Q$  is not a square  $((\alpha)\beta^{-1})^{(p-1)/2} \equiv  (\alpha^{p-1}/Q^{(p-1)/2}) \equiv -1$.
This phenomenon can be observed empirically by counting the proportion of $n$ which are lpsp for different options $P,Q$.
Those with $Q$ of the above form tend to be lpsp more than on average.

\begin{lemma}\label{lem4}
The period of the $V_{y+1}$ sequence divides the LCM of the orders
of $\alpha$ and $\beta$.
The order of $Q$ divides the LCM of the orders of $\alpha$ and $\beta$,
as does the order of $(\alpha)\beta^{-1}$.
All of these orders divide $p-1$.
\end{lemma}

\noindent
{\em Proof}.
It is clear from the definition of $V_{y+1}$ that its period divides
the periods of the two functions added to obtain it.
Since $Q=\alpha\beta$, its order must divide
the LCM of the orders of $\alpha$ and $\beta$.
Likewise, the order of $(\alpha)\beta^{-1}$ divides this LCM.
This completes the proof.

In contrast to the previous cases, if $Q\not\equiv1$ or $P/2\qmod{p}$,
there is no guarantee that there are any solutions at all.

However, we know that
$$V_{(p-2)+\ell(p-1)+1}\equiv2\qmod{p}
\quad\mathrm{and}\quad 
V_{\ell(p-1)+1}\equiv P\qmod{p}.$$

\begin{lemma}
 Let $t$ denote the order of $\alpha/\beta$ in $\mathbb{F}_p^\times$ then
 \[ V_{y+mt} = (\beta^t)^m V_y \]
 and hence each $a\in \mathbb{F}_p^\times$ is repeated by $V_y$ equally often as $(\beta^t)^ma$ and hence not more than $\ell/\ord(\beta^t)$ times within one period of $V_y\qmod{p}$.
 
 Let $f$ denote the order of $\alpha$ in $\mathbb{F}_p^\times$ then
 \[ V_{y+mf} = V_y + \beta^y((\beta^f)^m - 1) \]
 and hence each $a\in \mathbb{F}_p^\times$ is repeated by $V_y$ equally often as $(\beta^f)^ma$ and hence not more than $\ell/\ord(\beta^t)$ times within one period of $V_y\qmod{p}$.
 
 By symmetry each $a\in \mathbb{F}_p^\times$ is repeated by not more than $\ell/\ord(\alpha^{\ord(\beta)})$ times within one period of $V_y\qmod{p}$.
\end{lemma}

The proof is as in the previous case.

Now we consider the map
\[ \Psi : \mathbb{Z} \rightarrow {\mathbb F}_p \times {\mathbb F}_{p}^{\times} \]
given by
\[ y \mapsto  (V_y, N_{{\mathbb F}_{p^2}^{\times}/{\mathbb F}_{p}^{\times}}(\alpha^y)) =  (V_y,Q^y). \]
And note that this map has a period which divides the LCM of $\ord(\alpha)$ and $\ord(\beta)$, and equals it up to a multiple of $2$.

\begin{lemma}
 The function $\Psi$ is either $1:1$ or $2:1$ on its image, and hence each $a\in \mathbb{F}_p^\times$ is repeated by $V_y$ no more than $2\ord(Q)$ times within for $0\leq y \leq {\rm LCM}({\ord(\alpha)},{\ord(\beta)})$.
 
 It is $1:1$ unless there exists $k$ with $\beta=\alpha^k$ and $\alpha=\beta^k$.
\end{lemma}

The proof is as in the previous case.

\begin{proposition}
 The probability that $V_{py+1}\equiv 2Q \qmod{p}$ is less than the minimum of 
\[ \frac{2\ord(Q)}{{\rm LCM}({\ord(\alpha)},{\ord(\beta)})},\quad \frac{\ord(\alpha)}{{\rm LCM}({\ord(\alpha)},{\ord(\beta)})}, \frac{\ord(\beta)}{{\rm LCM}({\ord(\alpha)},{\ord(\beta)})}, \frac{\ord(\alpha/\beta)}{{\rm LCM}({\ord(\alpha)},{\ord(\beta)})}. \]
\end{proposition}

In contrast to the previous case it is challenging to get strong bounds on this expectation when the orders of $\alpha$ and $\beta$ are both large.
However, in that case one still expects the values of $\alpha^y$ and $\beta^y$ to behave like uniform random variables, and hence $V_y = \alpha^y + \beta^y$ should as well.

Note also that, as in the case $(D/p)=-1$ if we were considering the conditional probability of Congruence \eqref{eq2} given \eqref{eq1},
we would restrict to $y=(p-2)+x\ell$, where $\ell$ is the
order of $(\alpha)\beta^{-1}$ so that $U_{y+1}=0$, then as
$$(\alpha-\beta)U_z+2\beta^z=V_z$$
the condition
$$2Q=V_{y+1}=V_{(p-2)+x\ell+1}= 2\beta^{x\ell}$$
becomes
\[ Q = \beta^{x\ell}  \]
which is periodic with period the order of $\beta^\ell$, and is up to a multiple of $2$ the order of $\alpha\beta=Q$.
By symmetry we also obtain
\[ Q = \alpha^{x\ell}  \]
and taking the product of these two congruences gives
\[ Q^2 = Q^{x\ell} \]
and so $x\ell = 2 \qmod{{\rm ord}(Q)}$.
This final condition has a low probability of being satisfied if ${\rm ord}(Q)$ is large.
And is likely impossible to satisfy if $\ell = {\rm ord}((\alpha)\beta^{-1})$ is not relatively prime to ${\rm ord}(Q)$.

Once again we can see why the $V$ sequence probably
outperforms the $U$ sequence:
the expected period of the $U$ sequence is strictly less than $p-1$,
while that of $V$ is likely to be closer to $p-1$ as it comes from
the LCM of two periods.
The probability that $V_{y+1}=2Q$ can however be less than
$1/(p-1)$ since it can be 0.
Realistically, one expects that the values of the $V$ sequence
occur with equal frequency, though this is not guaranteed.

In contrast to the previous case the probabilities that for fixed $D$
and $p$ we have $U_{y+1}=0$ or $V_{y+1}=2Q$ are largely not
independent as they both depend on periods modulo $p-1$
and so likely either completely conflict or completely overlap.
This is consistent with the observation that there should be fewer
composites satisfying \eqref{eq1} and \eqref{eq2}
with $(D/n)=-1$ where most $p\mid n$ have $(D/p)=1$.

\subsection{Both cases together}\label{both}

We have seen that heuristically, if you
use \eqref{eq2} with $(D/p) = -1$ and allow
$Q = 1$ or $-1$ as in Method A, then the order of $\alpha$
divides $2(p+1)$.  But if you force $|Q| > 1$ as in Method A*,
then this order divides $p^2-1 = (p-1)(p+1)$ and probably
not $2(p+1)$ or $p-1$.

Consider the BPSW probable prime test.
The Fermat condition $2^{n-1} = 1 \qmod{n}$ basically
requires that the order of 2 modulo any prime factor
$p$ of $n$ divide $n-1$.  If the order is large, as it
often is, then it rarely divides $n-1$ when $n$ has
other prime factors than $p$.  But it is not that rare;
it does happen occasionally, and we get (some) pseudoprimes to base 2.

The Lucas condition \eqref{eq1} is trickier.  With $(D/n) = 1$,
it is just a Fermat test with 2 replaced by $\alpha$.
With $(D/n) = -1$, it operates in ${\mathbb F}_{p^2}$ and needs
$n+1$ to satisfy a congruence condition modulo $p+1$ in order to report $n$ probably prime.
Such a condition modulo $p+1$ happens about as often
as $p-1$ divides $n-1$, so we get occasional Lucas psps.

When we combine the Fermat and Lucas conditions, we
ask for $n$ to satisfy congruence conditions modulo both $p-1$ and $p+1$.
These probabilities are not independent, but as a first
approximation, each event has probability about $1/p$.
The probability of both events simultaneously would
be $1/p^2$.  Now $\sum_p 1/p$ diverges, while $\sum_p 1/p^2$
converges.  By the Borel-Cantelli lemmas, the first
event (just pseudoprimes or just Lucas pseudoprimes) occurs infinitely
often, while the second event (counterexample to BPSW)
occurs only finitely often.
(The second Borel-Cantelli Lemma requires the events be pairwise independent; the first
Lemma does not have this hypothesis.)

If we consider just \eqref{eq2} with $(D/n) = -1$ and use
a method for choosing $D$, $P$, $Q$ that does not allow
$|Q| = 1$, then we are forcing both congruence conditions on $n$ modulo $p-1$ and $p+1$  (or at least a large divisor
of $p-1$ respectively $p+1$), so the number of solutions should
be finite.  Of course, the fact that some orders
will be proper divisors of $p-1$, $p+1$, or $p^2-1$
will allow some solutions, perhaps infinitely
many, because the event probabilities are
greater than $1/p^2$ but still less than $1/p$.

\section{Are there infinitely many counterexamples to the test?} \label{S:counterexamples}

As noted earlier, no odd, composite $n$ is known that is psp(2) and, when Method A* is used to choose $P$ and $Q$, such that $n$ is also lpsp$(P, Q)$.

Nevertheless, if we search hard enough, we \emph{can} find odd, composite $n$ and Lucas parameters $P$ and $Q$ for which
\begin{align*}
2^{n-1} & \equiv 1 \qmod{n} \thinspace , \\
U_{n+1} & \equiv 0 \qmod{n} \thinspace ,  \quad \text{and} \\
V_{n+1} & \equiv 2Q \qmod{n}
\end{align*}
are simultaneously true.
One such example is $n = 341$, $P = 27$, $Q = 47$, $D = P^2 - 4Q = 541$.
This example was found by testing \emph{all possible} $(P, Q)$ pairs $\qmod{n}$, something one would not do when testing $n$ for primality.

We can also find odd, composite $n$ along with $P$ and $Q$, such that $n$ is simultaneously \emph{strong} psp(2), \emph{strong} lpsp$(P, Q)$, and vpsp$(P, Q)$.
This theorem was found empirically by testing all $(P, Q)$ pairs $\qmod{n}$.

\begin{theorem} \label{T:SpspSlpspVpsp}
Let $n \equiv 3 \qmod{4}$ be a strong pseudoprime base 2.
Let $k \ge 0$ be an integer.
Set $P = 2^k$ and $Q = 2^{2k-1}$.
Then $n$ is also a strong lpsp$(P, Q)$ and a vpsp$(P, Q)$.
\end{theorem}

Remarks.

1. If $k = 0$, then $Q = 2^{-1} \equiv (n+1)/2 \qmod{n}$.

2. Examples of spsp(2) that are $\equiv 3 \qmod{4}$ include \emph{composite} Mersenne numbers of the form $2^p - 1$ where $p$ is an odd prime \cite[p. 1008]{PSW1980}.

3. Corollary 1 below shows that infinitely many $n$ satisfy Theorem 1.

4. This proof uses the facts that $n$ is an epsp(2) and that $n \equiv 3 \qmod{4}$.
However, if $n \equiv 3 \qmod{4}$, then $n$ is epsp($a$) if and only if $n$ is spsp($a$) \cite[Theorem 4, p. 1009]{PSW1980}.

\medskip

\noindent
{\em Proof}.
First, $D = P^2 - 4Q = 2^{2k} - 4 \cdot 2^{2k-1} = 2^{2k} - 2 \cdot 2^{2k} = -(2^k)^2$.
Because $n \equiv 3 \qmod{4}$, this $D$ has Jacobi symbol $(D/n) = (-1/n) \cdot ((2^k)^2/n) = -1$.

Write $n+1 = d \cdot 2^s$, where $d$ is odd.
Then $s > 1$, and $2d \le (n + 1)/2$.
We will first prove that $V_{2d} \equiv 0 \qmod{n}$; by Congruence \eqref{E:StrongLrp2}, this will prove that $n$ is slpsp$(P, Q)$.

Let $\alpha$ and $\beta$ be the roots of the characteristic equation $x^2 - Px + Q = 0$, so that
\begin{gather*}
\alpha = \frac{P + \sqrt{D}}{2} = \frac{2^k + \sqrt{-4^k}}{2} = 2^{k-1} ( 1 + i) \thinspace , \\
\beta  = \frac{P - \sqrt{D}}{2} = \frac{2^k - \sqrt{-4^k}}{2} = 2^{k-1} ( 1 - i) \thinspace .
\end{gather*}

Then $V_{2d} = \alpha^{2d} + \beta^{2d}$.
Because $(1 + i)^2 = 2i$, we have
\[
\alpha^{2d}
 = \left( 2^{k-1} \right )^{2d} \cdot (1 + i)^{2d}
 = \left( 2^{2k-2} \right )^d \cdot (2i)^d
 = \left( 2^{2k-1} \right )^d \cdot i^d \thinspace .
\]
Similarly, because $(1 - i)^2 = -2i$, we have $\beta^{2d} = \left( 2^{2k-1} \right )^d \cdot (-i)^d$.
Therefore,
\[
V_{2d} = \alpha^{2d} + \beta^{2d} = \left( 2^{2k-1} \right )^d \cdot (i^d + (-i)^d) = 0 \thinspace ,
\]
so $n$ is a strong lpsp$(P, Q)$.

We will now prove that $V_{n+1} \equiv 2Q \qmod{n}$.
Because $n \equiv 3 \qmod{4}$, we can write $n + 1 = 4 M$, where $M$ is an integer.
Also, $V_{n+1} = \alpha^{n+1} + \beta^{n+1}$.

Observe that $(1 + i)^4 = (1 - i)^4 = -4$.
Then
\[
\alpha^{n+1}
= \left( 2^{k-1} \right )^{n+1} \cdot (1 + i)^{4M}
= \left( 2^{n+1} \right )^{k-1} \cdot (1 + i)^{4M}
= \left( 4 \cdot 2^{n-1} \right )^{k-1} \cdot (-4)^M \thinspace .
\]
$\beta^{n+1}$ has the same value.
Therefore,
\begin{align} \label{E:Vnp1}
V_{n+1} = \alpha^{n+1} + \beta^{n+1}
 & = 2 \cdot \left( 4 \cdot 2^{n-1} \right )^{k-1} \cdot (-1)^M \cdot 4^M \notag \\
 & = 2 \cdot 2^{2k-2} \cdot \left( 2^{n-1} \right)^{k-1} \cdot (-1)^M \cdot 2^{2M} \notag \\
 & = 2^{2k-1} \cdot \left( 2^{n-1} \right)^{k-1} \cdot (-1)^M \cdot 2^{(n+1)/2} \notag \\
 & = 2 Q \cdot \left( 2^{n-1} \right)^{k-1} \cdot (-1)^M \cdot 2^{(n-1)/2} \thinspace .
\end{align}
But $2^{n-1} \equiv 1 \qmod{n}$ because $n$ is spsp(2) and is therefore a Fermat pseudoprime base 2 that satisfies Congruence \eqref{E:FermatTest}.

Moreover, because $n$ is spsp(2), it is therefore an Euler pseudoprime base 2, so that, by Congruence \eqref{E:EulerCriterion},
$2^{(n-1)/2} \equiv \left( \frac{2}{n} \right) \qmod{n}$.
We now separate \eqref{E:Vnp1} into two cases.

Case I. $n \equiv 3 \qmod{8}$. Then (a), $M$ is odd, so $(-1)^M = -1$, and (b), $(2/n) = -1$.

Case II. $n \equiv 7 \qmod{8}$. Then (a), $M$ is even, so $(-1)^M = 1$, and (b), $(2/n) = 1$.

In both cases, \eqref{E:Vnp1} becomes
\[
V_{n+1} = 2 Q \cdot \left( 2^{n-1} \right)^{k-1} \cdot (-1)^M \cdot 2^{(n-1)/2}
   \equiv 2 Q \cdot 1 \cdot (-1)^M \cdot \left( \frac{2}{n} \right)
   \equiv 2 Q \qmod{n} \thinspace .
\]

Therefore, $n$ is also vpsp$(P, Q)$.
This completes the proof of the theorem.

With this $n$ and $Q$, the condition $Q^{(n+1)/2}\equiv Q \cdot (Q/n) \qmod{n}$
in step 5 of the enhanced primality test is also satisfied:
Since $n$ is a spsp(2), it is also an spsp($2^{2k-1}$), that is, spsp($Q$).
Therefore $n$ is an Euler pseudoprime to base 2:
$Q^{(n-1)/2}\equiv (Q/n) \qmod{n}$.  Multiply by $Q$ to get the
condition in step 5.

Note that the values of $n$ in Theorem \ref{T:SpspSlpspVpsp} are not counterexamples to our primality test, because Method A* never chooses these values of $P$ and $Q$.

Before we discovered paper \cite{PR80} we tried to prove on our own
that there are infinitely many spsp(2) in the congruence class $3\qmod{4}$
and found the theorem below, which has independent interest
and which needs the following lemma.

\begin{lemma}\label{lemmaQR}
For every positive integer $r$ there exists an integer $a\equiv3\qmod{4}$
such that for every odd prime $p$, if $p\equiv a\qmod{4r}$,
then $r$ is a quadratic residue modulo $p$: $\left(\frac{r}{p}\right)=+1$.
\end{lemma}

\noindent
{\em Proof}.
Write $r=2^st$ with $t$ odd.
If $t\equiv1\qmod{4}$, let $a=1+2t$.
If $t\equiv3\qmod{4}$, let $a=4t-1$.
In either case, if $s$ is odd, add $4t$ to $a$.
It is easy to see that $a\equiv3\qmod{4}$ in all cases.
In the rest of the proof suppose that $p$ is an odd prime
and $p\equiv a\qmod{4r}$.  Note that this implies that $p\equiv3\qmod{4}$.

If $r$ is a power of 4, then $s$ is even, $t=1$, $a=3$ and 
$\left(\frac{r}{p}\right)=\left(\frac{1}{p}\right)=+1$.

If $r$ is twice a power of 4, then $s$ is odd, $t=1$, $a=7$, $8\mid4r$ and 
$\left(\frac{r}{p}\right)=\left(\frac{2}{p}\right)=+1$
by a supplement to the Law of Quadratic Reciprocity (LQR) that says that if
$p\equiv 7\qmod{8}$, then $\left(\frac{2}{p}\right)=+1$.

Now suppose $s$ is even and $t>1$.
If $t\equiv1\qmod{4}$, then by the LQR we have
\begin{equation*}
\left(\frac{r}{p}\right)=
\left(\frac{2^s t}{p}\right)=
\left(\frac{t}{p}\right)=
\left(\frac{p}{t}\right)=
\left(\frac{a}{t}\right)=
\left(\frac{1+2r}{t}\right)=
\left(\frac{1}{t}\right)=+1.
\end{equation*}
If $t\equiv3\qmod{4}$, then by the LQR we have
\begin{equation*}
\left(\frac{r}{p}\right)=
\left(\frac{2^s t}{p}\right)=
\left(\frac{t}{p}\right)=
-\left(\frac{p}{t}\right)=
-\left(\frac{4t-1}{t}\right)=
-\left(\frac{-1}{t}\right)=+1.
\end{equation*}

Finally, suppose $s$ is odd and $t>1$.  Then $8\mid4r$ and $a\equiv7\qmod{8}$.
If $t\equiv1\qmod{4}$, then by the LQR we have
\begin{equation*}
\left(\frac{r}{p}\right)=
\left(\frac{2^s t}{p}\right)=
\left(\frac{2}{p}\right) \left(\frac{t}{p}\right)=
(+1) \left(\frac{p}{t}\right)=
\left(\frac{a}{t}\right)=
\left(\frac{1+2t+4t}{t}\right)=
\left(\frac{1}{t}\right)=+1.
\end{equation*}
If $t\equiv3\qmod{4}$, then by the LQR we have
\begin{equation*}
\left(\frac{r}{p}\right)=
\left(\frac{2^s t}{p}\right)=
\left(\frac{2}{p}\right) \left(\frac{t}{p}\right)=
-\left(\frac{p}{t}\right)=
-\left(\frac{a}{t}\right)=
-\left(\frac{4t-1+4t}{t}\right)=
-\left(\frac{-1}{t}\right)=+1.
\end{equation*}

This completes the proof.

\medskip

\begin{theorem}\label{thmspspr}
If $r>1$ is an integer, there are infinitely many Carmichael numbers
$m\equiv3\qmod{4}$ that are also strong pseudoprimes to base $r$.
Moreover, there is a constant $K>0$ which depends on $r$
so that the number of such $m<X$
is $\ge X^{K/(\log\log\log X)^2}$ for all sufficiently large $X$.
\end{theorem}

\noindent
{\em Proof}.
Let $M=4r$ and choose $a$ by Lemma \ref{lemmaQR}.
Wright \cite{Wri13} proved that there are infinitely many
Carmichael numbers $m\equiv a\qmod{M}$ and in fact
$\ge X^{K/(\log\log\log X)^2}$ of them below $X$ for all large enough $X$.
Since $a\equiv3\qmod{4}$ and $4\mid M$, we have $m\equiv 3\qmod{4}$.
In the construction of the Carmichael numbers $m$ in the proof
in \cite{Wri13}, every prime factor $p$ of $m$ is odd and $\equiv a\qmod{M}$.
Thus each $p\mid m$ satisfies $p\equiv 3\qmod{4}$ and, by Lemma \ref{lemmaQR},
$\left(\frac{r}{p}\right)=+1$.
By Corollary 1.2 of \cite{AGP94}, if $\left(\frac{r}{p}\right)$
has the same value for every prime $p\mid m$, then $m$ is a strong
pseudoprime to base $r$.
This completes the proof.

\medskip

In 1980, van der Poorten and Rotkiewicz \cite{PR80} proved that for every
integer $r>1$ there are infinitely many spsp($r$) in every
arithmetic progression $ax+b$ with $(a,b)=1$, but they did not
bound the growth rate of such numbers.

Theorem 1 of \cite{PSW1980} asserts that for all $r>1$ and
$x>r^{15r}+1$, there are more than $(\log x)/(4r\log r)$
strong pseudoprimes to base $r$ less than $x$.
Every one of these spsp($r$) is $\equiv1\qmod{4}$.
Later, Pomerance \cite{Pom82} proved an even greater lower bound
on the number of strong pseudoprimes to base $r$ less than $x$.
All of the spsp($r$) he constructed are $\equiv1\qmod{4}$.

\medskip

\begin{corollary}\label{cor3mod4}
Let $k$ be a nonnegative integer.  Let $P=2^k$ and $Q=2^{2k-1}$.
Then there exist infinitely many Carmichael numbers $m\equiv3\qmod{4}$
that are strong pseudoprimes to base $2$, strong lpsp($P,Q$)
and vpsp($P,Q$).
Moreover, there is a constant $K>0$ so that the number of such $m<X$
is $\ge X^{K/(\log\log\log X)^2}$ for all sufficiently large $X$.
\end{corollary}

\medskip

The corollary follows from Theorems \ref{T:SpspSlpspVpsp} and
\ref{thmspspr}.

In the case $r=2$ all spsp(2) that we constructed in Theorem
\ref{thmspspr} are $\equiv7\qmod{8}$.
This is because when $r=2$ Lemma \ref{lemmaQR} sets $s=t=1$ and $a=7$.
It is easy to modify the proof of Theorem \ref{thmspspr} to show
that there are infinitely many spsp(2) that are $\equiv3\qmod{8}$.
Rather than use Lemma \ref{lemmaQR}, just set $M=8$ and $a=3$.
Then Wright's proof for this arithmetic progression
constructs many Carmichael numbers $m\equiv3\qmod{8}$,
every prime factor of which is also $\equiv3\qmod{8}$.
Then the Legendre symbols $\left(\frac{2}{p}\right)$ are all $-1$
by a supplement to the LQR
so that Corollary 1.2 of \cite{AGP94} still applies to show
that $m$ is spsp(2).

The smallest Carmichael number that satisfies all the conditions of this Corollary is 3\,215\,031\,751.

\medskip

\textbf{Pomerance's heuristic argument.}

In 1899, Korselt \cite{Kor99} proved that $n$ is a Carmichael number if and only if
$n$ is square free, has at least three prime factors, and for each
prime $p$ dividing $n$ we have $p-1$ divides $n-1$.
As part of Erd{\H o}s’ heuristic argument (see \cite{Erd56}) that there are infinitely many Carmichael numbers (in fact more than $x^{1-\epsilon}$ up to $x$),
he shows that there are many composite square free numbers $n$ for which $p-1$ divides $n-1$ for each prime factor $p$ of $n$.
Pomerance \cite{Pom84} modified this argument to show that there are infinitely
many strong pseudoprimes $n$ to base 2 that are also Lucas pseudoprimes
with $(D/n)=-1$.
Pomerance's argument showed heuristically that there are Carmichael numbers $n$
so that for each prime $p$ dividing $n$ we have $p+1$ divides $n+1$,
and congruence conditions to ensure that $n$ is a strong pseudoprime
to base 2 and $(D/n)=(5/n)=-1$.
Since the numbers $n$ Pomerance constructed all satisfy $p-1\mid n-1$
and $p+1\mid n+1$ for each prime factor $p$ of $n$, they
satisfy all of Congruences \eqref{eq1}--\eqref{eq4}.

Pomerance chooses an integer $k>4$ and a large $T$.
He lets $P_k(T)$ be the set of all primes $p$ in $[T,T^k]$
such that
\begin{enumerate}
\item $p\equiv3\qmod{8}$ and the Jacobi symbol $(5/p)=-1$.
\item $(p-1)/2$ is square free and composed only of primes $q<T$ with
$q\equiv1\qmod{4}$.
\item $(p+1)/4$ is square free and composed only of primes $q<T$ with
$q\equiv3\qmod{4}$.
\end{enumerate}

Let $Q_1$ be the product of all primes $q<T$ with $q\equiv1\qmod{4}$.
Let $Q_3$ be the product of all primes $q<T$ with $q\equiv3\qmod{4}$.

Heuristically, the size of $P_k(T)$ is about $T^k/\log^2 T$.

Let $\ell$ be odd and let $n$ be any product of $\ell$ primes $p\in P_k(T)$
such that $n\equiv1\qmod{Q_1}$ and $n\equiv-1\qmod{Q_3}$.

Then $n\equiv3\qmod{8}$, $(5/n)=-1$ and for all primes $p\mid n$
we have $p-1\mid n-1$ and $p+1\mid n+1$.
This implies that $n$ is a strong pseudoprime to base 2
and $n$ satisfies all of \eqref{eq1}--\eqref{eq4},
so $n$ is a Lucas pseudoprime, a v pseudoprime, and so is a counterexample to the enhanced BPSW primality test.

The arguments of both Erd\H{o}s and Pomerance were heuristic,
with many unproved but plausible assumptions.

The condition $k>4$ allows one to show that there are $x^{1-\epsilon}$
counterexamples $n$ to the enhanced BPSW primality test with $n<x$.

The conditions $p\equiv3\qmod{8}$ for $p\in P_k(T)$
make it easy to prove $n$ is spsp(2).

A computer search for counterexamples to BPSW using Pomerance's
construction would be very slow due partly to the conditions
2 and 3 above for primes $q<T$ and due partly to the
conditions $n\equiv1\qmod{Q_1}$ and $n\equiv-1\qmod{Q_3}$.

\medskip

{\bf Conclusion.}

In Section \ref{S:Qgt1}, we have presented some reasons why
counterexamples to the BPSW test or to our new strengthened
test should be rare or nonexistent.  
On the other hand, in this section we have suggested
that there might be many, perhaps infinitely many
counterexamples to these tests.
So which is it?

The arguments in Section \ref{S:Qgt1} seem to apply to relatively
small numbers, those with hundreds or thousands of decimal
digits that we might actually test for primality using
computers.
We believe that counterexamples to either test are extremely
rare among numbers of that size.
The arguments in this section seem to apply to truly enormous
numbers, numbers too large for even a computer to multiply.
Some day, when we know more about the distribution of primes,
we might be able to prove rigorously that there are infinitely
many counterexamples, but these numbers might be so large that
their logarithm exceeds the number of electrons in the universe.

\section{Open Questions} \label{S:OpenQuestions}

Suppose $n$ is composite and that we have the full factorization of $n$.
Work by various authors has produced formulas that count or estimate:

\begin{itemize}
\item the number of bases $a$ for which $n$ is a psp($a$) \cite[Thm. 1]{{BW80}}, \cite[Lemma 1]{Monier1980}
\item the number of bases $a$ for which $n$ is a spsp($a$) \cite[Prop. 1]{Monier1980}, \cite[Thm. 1.4]{Arn97}
\item given $D$, the number of $P$ for which there is a $Q$ such that $n$ is lpsp$(P, Q)$ \cite[Thm. 2]{{BW80}}
\item given $D$, the number of $(P, Q)$ pairs for which $n$ is slpsp$(P, Q)$ \cite[Thm. 1.5]{Arn97}
\item the number of $(P, Q)$ pairs $\qmod{n}$ for which $n$ is simultaneously lpsp$(P, Q)$ and vpsp$(P, Q)$ \cite[Thm. 16]{FioriShallue}.
\end{itemize}

We would like to have a formula that bounds, or better yet, counts, the number of $D$, or the number of $(P, Q)$ pairs for which $n$ is a vpsp$(P, Q)$.
We would also like to see an estimate of the asymptotic growth rate for the number of vpsp's $\leq x$; this would presumably depend on the algorithm for choosing $P$ and $Q$.

\bigskip

\appendix  
\section{Appendix. Methods A and A* generate the same lpsp lists} \label{S:Astar}

Recall that a Lucas probable prime is a solution $n$ to \eqref{E:UCongruence}.
Here we prove that Methods A and A* give the same solutions.
We also prove a similar result for strong Lucas probable primes.
In this appendix we write $U_n(P,Q)$ and $V_n(P,Q)$ for the
two Lucas sequences with parameters $P$, $Q$.

\begin{theorem}
Let $n$ be a positive integer relatively prime to 10.
Then $n$ is a Lucas probable prime for Method A if and only if
it is a Lucas probable prime for Method A*.
\end{theorem}

\noindent
{\em Proof}.
Methods A and A* differ only when $D=5$.
With that $D$, Method A sets $P=1$, $Q=-1$, while A* sets $P=Q=5$.
Let
$$\alpha_1=\frac{1+\sqrt{5}}{2}\quad\mathrm{and}\quad 
\beta_1=\frac{1-\sqrt{5}}{2}$$
be the two roots of $x^2-x-1=0$ and let
$$\alpha_2=\frac{5+\sqrt{5}}{2}\quad\mathrm{and}\quad 
\beta_2=\frac{5-\sqrt{5}}{2}$$
be the two roots of $x^2-5x+5=0$.
Then
\begin{equation}
\alpha_2^2=5\left(\frac{3+\sqrt{5}}{2}\right)=5\alpha_1^2
\quad\mathrm{and}\quad 
\beta_2^2=5\left(\frac{3-\sqrt{5}}{2}\right)=5\beta_1^2. \label{eqab}
\end{equation}
Since $\alpha_1-\beta_1=\sqrt{5}=\alpha_2-\beta_2$, we have
$$U_{2k}(5,5)=\frac{\alpha_2^{2k}-\beta_2^{2k}}{\alpha_2-\beta_2}=
5^k\left(\frac{\alpha_1^{2k}-\beta_1^{2k}}{\alpha_1-\beta_1}\right)=
5^kU_{2k}(1,-1).$$
Now $n$ is odd, so we can write $n + 1 = 2k$; by the previous expression, $$U_{n+1}(5,5) = 5^{(n+1)/2}U_{n+1}(1,-1).$$
Since $n$ is not a multiple of 5, $U_{n+1}(5,5)\equiv0\qmod{n}$
if and only if $U_{n+1}(1,-1)\equiv0\qmod{n}$.
This completes the proof.

Now we prove the analogue of this theorem for the strong lprp test.

\begin{theorem}
Let $n$ be a positive integer relatively prime to 10.
Then $n$ is a strong Lucas probable prime for Method A if and only if
it is a strong Lucas probable prime for Method A*.
\end{theorem}

\noindent
{\em Proof}.
Let $\alpha_1$, $\alpha_2$, $\beta_1$, $\beta_2$ be as in the
proof of the previous theorem.
We will prove that
\begin{align}
U_{2k+1}(5,5) & = 5^kV_{2k+1}(1,-1), \label{eq11} \\
V_{2k+1}(5,5) & = 5^{k+1}U_{2k+1}(1,-1),\quad\mathrm{and}\label{eq12} \\
V_{2k}(5,5) & = 5^kV_{2k}(1,-1), \label{eq13}
\end{align}
Using Equation \eqref{eqab},
the left side of \eqref{eq11} is
$$\frac{\alpha_2^{2k+1}-\beta_2^{2k+1}}{\alpha_2-\beta_2} =
\frac{5^k}{\alpha_2-\beta_2}\left(\alpha_1^{2k}\alpha_2 -
\beta_1^{2k}\beta_2\right) =
\frac{5^k}{\sqrt{5}}\left(\alpha_1^{2k+1}\frac{\alpha_2}{\alpha1} -
\beta_1^{2k+1}\frac{\beta_2}{\beta_1}\right).$$
Now
$$\frac{\alpha_2}{\alpha_1} = \frac{5+\sqrt{5}}{1+\sqrt{5}}=\sqrt{5}
\quad\mathrm{and}\quad 
\frac{\beta_2}{\beta_1} = \frac{5-\sqrt{5}}{1-\sqrt{5}}=-\sqrt{5},$$
so the left side of \eqref{eq11} becomes
$$\frac{5^k}{\sqrt{5}}\left(\alpha_1^{2k+1}\sqrt{5} -
\beta_1^{2k+1}(-\sqrt{5}) \right) =
5^k\left(\alpha_1^{2k+1} + \beta_1^{2k+1}\right),$$
which is the right side of \eqref{eq11}.
Equation \eqref{eq12} is proved the same way.
Equation \eqref{eq13} is even easier:
$$ V_{2k}(5,5) = \alpha_2^{2k} + \beta_2^{2k} =
(5\alpha_1)^{2k} + (5\beta_1)^{2k} = 5^k(\alpha_1^{2k} + \beta_1)^{2k} = 
5^kV_{2k}(1,-1).$$

Now if $n$ is an slprp(1,$-1$) because $V_d(1,-1) \equiv 0 \qmod{n}$,
then Equation \eqref{eq11} shows that $n$ is an slprp(5,5) because
$U_d(5,5) \equiv 0 \qmod{n}$, and vice versa.
Also if $n$ is an slprp(1,$-1$) because $U_d(1,-1) \equiv 0 \qmod{n}$,
then Equation \eqref{eq12} shows that $n$ is an slprp(5,5) because
$V_d(5,5) \equiv 0 \qmod{n}$, and vice versa.
Finally, if $n$ is an slprp(1,$-1$) because $V_{d \cdot 2^r}(1,-1) \equiv 0 \qmod{n}$
for some $0 < r < s$,
then Equation \eqref{eq13} shows that $n$ is an slprp(5,5) because
$V_{d \cdot 2^r}(5,5) \equiv 0 \qmod{n}$, and vice versa.
This completes the proof.

\end{document}